\documentclass[12pt]{amsart}

\usepackage{amscd,amsmath,amssymb}
\usepackage[mathscr]{eucal}

\makeatletter \@addtoreset{equation}{section}\makeatother

\newtheorem{theorem}{Theorem}[section]
\newtheorem{lemma}[theorem]{Lemma}
\newtheorem{proposition}[theorem]{Proposition}

\title{\bf On Serre duality for compact homologically smooth DG algebras}
\author{D.Shklyarov}

\begin{document}

\maketitle

\centerline{\it{\small To Leonid L'vovich Vaksman on his 55th
birthday, with gratitude}}

\baselineskip 1.50pc

\section{Introduction}

Let $X$ be a smooth projective variety over a perfect field $k$. It
is a classical fact that the bounded derived category
$\mathsf{D}^b(\mathsf{Coh}X)$ of the category of coherent sheaves on
$X$ is equivalent (as a triangulated category) to the derived
category $\mathsf{D}_{\mathsf{per}}(A)$ of perfect modules over a DG
(=differential graded) algebra $A$ (see \cite{BVDB, Ro} and
references therein). The equivalence
$\mathsf{D}^b(\mathsf{Coh}X)\simeq\mathsf{D}_{\mathsf{per}}(A)$
implies some properties of $A$. First of all, the total cohomology
of $A$ has to be finite-dimensional (such DG algebras are called
compact \cite{G1,KS}). Smoothness of $X$ boils down to $A$ being
perfect as an $A$-bimodule\footnote{It seems that a rigorous proof
of this statement is yet to be written.} (such DG algebras are
called homologically smooth \cite{G1,KS}; see also \cite{TV}).

In view of the above observation, it is natural to expect that
compact homologically smooth DG algebras possess many properties of
smooth projective varieties. For example, if some invariant of
smooth projective varieties depends on the derived category
$\mathsf{D}^b(\mathsf{Coh}X)$, rather then on $X$ itself, then one
can try to define and study the corresponding invariant for compact
homologically smooth DG algebras. The aim of this approach is
two-fold. On one hand, one can try to obtain `geometry-free` proofs
of some classical results which may prove useful for the purposes of
noncommutative algebraic geometry. On the other hand, some of the
classical results (e.g. collapsing of the Hodge-to-de Rham spectral
sequence), if true in this noncommutative setting, would have some
interesting and important applications (e.g. to topological
conformal field theories \cite{C,KS}).

In this paper, we deal with an algebraic counterpart of the Hodge
cohomology groups $\mathrm{H}^{p}(X, \Omega^q_X)$ which is the
Hochschild homology $\mathrm{HH}_n(A)$ of a DG algebra $A$.

The Hochschild homology was originally defined for ordinary
associative algebras. It has been generalized in different
directions: there is a definition of the Hochschild homology of any
DG algebra or any scheme which agree with that of ordinary algebras.
It is a consequence of deep results of \cite{K1,K11/2,K2} that the
equivalence
$\mathsf{D}^b(\mathsf{Coh}X)\simeq\mathsf{D}_{\mathsf{per}}(A)$
implies an isomorphism of the Hochschild homology groups
$\mathrm{HH}_n(X)=\mathrm{HH}_n(A)$. It turns out that for a smooth
projective variety $X$, one has an isomorphism
$\mathrm{HH}_n(X)=\oplus_i \mathrm{H}^{i-n}(X, \Omega^i_X)$. Thus,
$\mathrm{HH}_n(A)$ is a 'right' replacement of the Hodge cohomology
in the general setting.

The first fundamental property of the Hodge cohomology of smooth
projective varieties is finite-dimensionality:
$\sum_{p,q}\mathrm{dim}\,\mathrm{H}^{p}(X, \Omega^q_X)<\infty$.
Therefore, it is natural to expect that
$\sum_n\mathrm{dim}\,\mathrm{HH}_n(A)<\infty$ for an arbitrary
compact homologically smooth DG algebra $A$. Furthermore, by the
classical Serre duality, there exists a non-degenerate pairing
$\mathrm{H}^{p}(X, \Omega^q_X)\times\mathrm{H}^{s}(X,
\Omega^{t}_X)\rightarrow k$ provided $p+s=q+t=\dim\,X$. Again, we
may hope that there exists a non-degenerate pairing
$\mathrm{HH}_n(A)\times \mathrm{HH}_{-n}(A)\rightarrow k$ on the
algebraic side. The aim of the paper is to prove the above two
assertions (see Theorem \ref{main}). Both results, we believe, are
well-known to the experts (cf. \cite{KS}). For the case of
associative algebras, they can be derived from \cite{vdb}. A similar
"categorical" approach to Serre duality in the geometric setting can
be found in \cite{Ma0,Ma,Ca}.

We would like to point out one interesting corollary of the
existence of the pairing on $\mathrm{HH}_\bullet(A)$. Namely, if $A$
is a compact homologically smooth DG algebra concentrated in {\it
non-positive} degrees then its Hochschild homology is concentrated
in degree 0. Indeed, it follows from the explicit form of the
bar-resolution of a DG algebra $A$ (see, for example, \cite{GJ})
that $\mathrm{HH}_n(A)=0$ for $n>0$ provided $A$ is concentrated in
non-positive degrees. This, together with the existence of a
non-degenerate pairing $\mathrm{HH}_n(A)\times
\mathrm{HH}_{-n}(A)\rightarrow k$, implies the result. This
corollary gives, for example, an alternative proof of the main
result of \cite{Ci} which describes the Hochschild homology of some
quiver algebras with relations.

There is yet another application of the corollary. It is related to
an analog of the afore-mentioned collapsing of the Hodge-to-de Rham
spectral sequence. By the classical Hodge theory, the de Rham
differential vanishes on the Hodge cohomology of a smooth projective
variety. An algebraic counterpart of this fact is the following
statement: Connes' differential
$B:\mathrm{HH}_\bullet(A)\to\mathrm{HH}_{\bullet-1}(A)$ vanishes
whenever $A$ is compact and homologically smooth. This statement is
the well known Noncommutative Hodge-to-de Rham degeneration
conjecture formulated by M. Kontsevich and Y. Soibelman several
years ago (actually, the conjecture is a stronger statement than
just vanishing of $B$ on the Hochschild homology; its precise
formulation can be found in \cite{KS}). Recently, D. Kaledin
\cite{Ka} proved this conjecture in the special (although very
difficult) case of DG algebras concentrated in non-negative degrees.
We would like to notice that the above corollary of our main result
implies the conjecture in the much easier case of homologically
smooth compact DG algebras concentrated in non-positive degrees.
Indeed, the Hochschild homology of such a DG algebra is concentrated
in degree 0 and $B$ is bound to vanish.

\pagebreak

\noindent{\bf Acknowledgments}

First and foremost, I wish to thank Yan Soibelman for introducing me
to this very interesting emerging part of Noncommutative Algebraic
Geometry. I am also very grateful to Dmitry Kaledin and Bernhard
Keller for taking the time to read a preliminary version of this
paper and making a number of helpful remarks. Of course, I am alone
responsible for any possible mistakes and inaccuracies.

\medskip
\noindent{\bf Notation}

Throughout the paper, we work over a fixed field $k$. All vector
spaces, algebras, linear categories are defined over $k$. All the
definitions regarding DG algebras and DG categories we are going to
use can be found in \cite{K}.

If $A$ is a DG algebra
$$
A=\bigoplus\limits_{n\in\mathbb{Z}}A^n,\quad d=d_A: A^n\rightarrow A^{n+1}
$$
then $A^\mathsf{ op}$ stands for the opposite DG algebra. We denote
the DG algebra $A^\mathsf{ op}\otimes A$ with the standard DG
algebra structure by $A^\mathsf{e}$.

The DG category of right DG $A$-modules is denoted by $\mathsf{
Mod}(A)$. We will write $n.a$ for the action of an element $a\in A$
on an element $n$ of a right DG $A$-module $N$. The homotopy and the
derived categories of $\mathsf{ Mod}(A)$ are denoted by $\mathsf{
Ho}(A)$ and $\mathsf{ D}(A)$, respectively.

Except for few cases, we work exclusively with {\it right} DG
modules over DG algebras. Therefore, we don't reserve a notation for
the category of left DG modules. We use the fact that right
$A^\mathsf{ op}$-modules are canonically left $A$-modules. For
example, we consider the tensor product $N{\otimes}_AM$ of a right
$A$-module $N$ and a right $A^\mathsf{ op}$-module $M$.

Let us recall the notation for some standard derived functors.  Let
$A$ and $B$ be two DG algebras. Take $N\in \mathsf{Mod}(A)$ and
$X\in \mathsf{Mod}(A\otimes B)$. The vector space
$\mathrm{Hom}_{\mathsf{Mod}(A)}(N, X)$ is canonically a right DG
$B$-module:
\begin{equation}\label{actia}
f.b(n)=(-1)^{|b||n|}f(n).b.
\end{equation}
The assignment $(N, X)\mapsto\mathrm{Hom}_{\mathsf{Mod}(A)}(N, X)$
is a DG bifunctor. It gives rise to a triangulated bifunctor
$\mathsf{Ho}(A)^\mathsf{op}\times \mathsf{Ho}(A\otimes B)\rightarrow
\mathsf{Ho}(B)$ denoted by $\mathsf{Hom}_{A}(-, -)$. Fix $X$. The
corresponding derived functor $\mathsf{D}(A)^\mathsf{op}\rightarrow
\mathsf{D}(B)$ is denoted by $\mathsf{RHom}_{A}(-, X)$. It is
defined by means of appropriate projective resolutions \cite{K}. If
$B=k$ then $\mathrm{H}^0(\mathsf{RHom}_{A}(N,
X))\simeq\mathrm{Hom}_{\mathsf{D}(A)}(N, X)$.

Let $C$ be yet another DG algebra. Given $N\in\mathsf{Mod}(C\otimes
A)$ and $M\in \mathsf{Mod}(A^\mathsf{op}\otimes B)$, consider
$N{\otimes}_AM\in \mathsf{Mod}(C\otimes B)$. The assignment $(N,
M)\rightarrow N{\otimes}_AM$ is a DG bifunctor. It induces a triangulated
bifunctor $\mathsf{Ho}(C\otimes A)\times
\mathsf{Ho}(A^\mathsf{op}\otimes B)\rightarrow \mathsf{Ho}(C\otimes
B)$ which we denote by $-\dot{\otimes}_A-$.  The total derived
functor $\mathsf{D}(C\otimes A)\times
\mathsf{D}(A^\mathsf{op}\otimes B)\rightarrow \mathsf{D}(C\otimes
B)$ is denoted by $-\otimes^\mathsf{L}_A-$. It is also defined by
means of projective resolutions.

\medskip

\section{Perfect modules and the notion of smoothness}\label{pm}
We start by recalling the construction of the DG category $\mathsf{
Tw}(A)$ of twisted modules over $A$ (see \cite{BK}).

Let us view $A$ as a DG category with one object. The first step is
to enlarge $A$ to a new DG category $\mathbb{Z}A$ by adding formal
shifts of the object. Namely, the objects of $\mathbb{Z}{A}$ are
enumerated by integers and denoted by $A[n]$, $n\in\mathbb{Z}$. The
space of morphisms $\mathrm{Hom}_{\mathbb{Z}{A}}(A[n], A[m])$
coincides with $\mathrm{Hom}_{A}(A, A)[m-n]\simeq A[m-n]$ as a
$\mathbb{Z}$-graded vector space; the differential
$d_{\mathbb{Z}{A}}$ on $\mathrm{Hom}_{\mathbb{Z}{A}}(A[n], A[m])$ is
given by $d_{\mathbb{Z}{A}}= (-1)^md_{A}$.

The objects of the category $\mathsf{ Tw}(A)$ are pairs
$(\bigoplus\limits_{j=1}^nA[r_j], \alpha)$, where $r_1, r_2, \ldots
r_n$ are integers, $\alpha=(\alpha_{ij})$ is a strictly upper
triangular $n\times n$-matrix of morphisms $\alpha_{ij}\in
\mathrm{Hom}_{\mathbb{Z}{A}}^1(A[r_j], A[r_i])\simeq A^{1+r_i-r_j}$
satisfying the Maurer-Cartan equation
$$
d_{\mathbb{Z}{A}}(\alpha)+\alpha\cdot\alpha=0.
$$
The space $\mathrm{Hom}_{\mathsf{
Tw}(A)}((\bigoplus\limits_{j=1}^nA[r_j], \alpha),
(\bigoplus\limits_{i=1}^mA[s_i], \beta))$ is a $\mathbb{Z}$-graded
space whose $p$-th component consists of $m\times n$-matrices
$f=(f_{ij})$ of morphisms $f_{ij}\in
\mathrm{Hom}_{\mathbb{Z}{A}}^p(A[r_j], A[s_i])=A^{p+s_i-r_j}$. The
differential $d_{\mathsf{ Tw}(A)}$  is defined as follows
$$
d_{\mathsf{ Tw}(A)}(f)=d_{\mathbb{Z}{A}}(f)+\beta\cdot f-(-1)^{|f|}f\cdot\alpha.
$$

The direct sum on objects of the category $\mathsf{ Tw}(A)$ is
defined in the obvious way:
$$
(\oplus A[r_j], \alpha)\bigoplus(\oplus A[s_i], \beta)=(\oplus
A[r_j]\oplus A[s_i], \begin{pmatrix}\alpha & 0 \\ 0 & \beta
\end{pmatrix}).
$$
To make $\mathsf{ Tw}(A)$ additive, we add also the zero object 0
satisfying the usual axioms. $A$ is embedded into  $\mathsf{ Tw}(A)$
as a full DG subcategory via the identification $A=(A, 0)$.

The homotopy category  $\mathsf{ Ho}(\mathsf{ Tw}(A))$ is equipped
with a canonical triangulated structure as follows. The shift
functor on $\mathsf{ Tw}(A)$ is given by
$$
(\oplus A[r_j], \alpha)[1]=(\oplus A[r_j+1], -\alpha)
$$
(the shift acts trivially on morphisms). Suppose
$f\in\mathrm{Hom}_{\mathsf{ Tw}(A)}^0((\oplus A[r_j], \alpha),
(\oplus A[s_i], \beta))$ is closed, i.e. $d_{\mathsf{ Tw}(A)}(f)=0$.
The cone of $f$, $\mathrm{Cone}(f)$, is the object of $\mathsf{
Tw}(A)$ defined by
$$
\mathrm{Cone}(f)=(\oplus A[s_i]\oplus A[r_j+1], \begin{pmatrix}\beta
& f \\ 0 & -\alpha \end{pmatrix}).
$$
By definition, a triangle in $\mathsf{ Ho}(\mathsf{ Tw}(A))$ is
distinguished if it is isomorphic to a triangle of the form
\[X\stackrel{f}{\longrightarrow}Y\stackrel{j}{\longrightarrow}\mathrm{Cone}(f)\stackrel{p}{\longrightarrow}X[1]\]
where $f$ is a degree 0 closed morphism and $j$ and $p$ are defined
in the obvious way.

The category $\mathsf{ Tw}(A)$ can be embedded into the DG category
$\mathsf{ Mod}(A)$ by means of the so called Yoneda embedding
\cite{BK}. Let us describe the image of  $(\bigoplus A[r_j],
\alpha)$ under $Y$. We can replace the formal direct sum by the
usual direct sum and view $X=\bigoplus A[r_j]$ as a usual graded
space. The matrix $\alpha$ allows us to equip $X$ with the structure
of a DG module by $d_X=d+\alpha.$ Here $d$ acts on $A[r_j]$ as the
differential of $A$ multiplied by $(-1)^{r_j}$. Finally, the right
$A$-module structure is given by the coordinate-wise multiplication
from the right:
$$
[a_1, a_2,\ldots].a=[a_1a, a_2a,\ldots], \quad a\in A,\quad a_j\in A[r_j].
$$

The functor $Y$ induces an embedding  $\mathsf{Ho}(\mathsf{
Tw}(A))\rightarrow \mathsf{Ho}(A)$ which allows us to view
$\mathsf{Ho}(\mathsf{ Tw}(A))$ as a triangulated subcategory of
$\mathsf{Ho}(A)$. Let $\mathsf{Ho}_{\mathsf{per}}(A)$ be the
smallest full triangulated subcategory of $\mathsf{Ho}(A)$
containing $\mathsf{Ho}(\mathsf{Tw}(A))$ and closed under taking
isomorphisms and direct summands. It embeds into a larger
triangulated subcategory
$\overline{\mathsf{Ho}}_{\mathsf{per}}(A)\subset\mathsf{Ho}(A)$
defined as follows: $N\in\overline{\mathsf{Ho}}_{\mathsf{per}}(A)$
if there exists a distinguished triangle in $\mathsf{Ho}(A)$ of the
form
\[pN\longrightarrow N\longrightarrow aN\longrightarrow pN[1]\]
where $pN\in\mathsf{Ho}_{\mathsf{per}}(A)$ and $aN$ is acyclic. In
other words, if we denote the triangulated subcategory of acyclic
modules by $\mathsf{Ho}_{\mathsf{ac}}(A)$ then
$$
\overline{\mathsf{Ho}}_{\mathsf{per}}(A)=\mathsf{Ho}_{\mathsf{per}}(A)\star\mathsf{Ho}_{\mathsf{ac}}(A)
$$
(see \cite{BVDB} for a description of the operation $\star$). We
define the triangulated subcategory
$\mathsf{D}_{\mathsf{per}}(A)\subset\mathsf{D}(A)$ of perfect
modules to be the image of
$\overline{\mathsf{Ho}}_{\mathsf{per}}(A)$ under the localization
$\mathsf{Ho}(A)\rightarrow\mathsf{D}(A)$.

The above chain of definitions can be rephrased as follows: a
perfect DG $A$-module is a module which can be resolved by a direct
summand of a twisted module.

Let us recall several useful results about perfect modules which we
are going to use in the sequel.

A DG $A$-module $P$ is called homotopically projective if
$\mathrm{Hom}_{\mathsf{Ho}(A)}(P, N)=0$ for any acyclic module $N$
(such modules are also called $K$-projective \cite{BL}). It is known
that objects of $\mathsf{Ho}_{\mathsf{per}}(A)$ are homotopically
projective modules \cite[\S 13]{D}. This observation allows us to
apply some well known facts about homotopically projective modules
to perfect modules. For example, one has \cite[Corollary
10.12.4.4]{BL}:
\begin{proposition}\label{flat}
If $P$ is homotopically projective then $P\dot{\otimes}_AN$ is
acyclic whenever $N$ is an acyclic left DG $A$-module.
\end{proposition}

Here is one more application. It is a classical fact that the
restriction of the localization
$\mathsf{Ho}(A)\rightarrow\mathsf{D}(A)$ onto the homotopy
subcategory of homotopically projective modules is an equivalence
(the quasi-inverse functor sends modules to their homotopically
projective resolutions).  Therefore $\mathsf{D}_{\mathsf{per}}(A)$
is equivalent to $\mathsf{Ho}_{\mathsf{per}}(A)$.

Finally, we are ready to formulate the main definition of this
paper: A DG algebra $A$ is said to be {\bf homologically smooth} if
$A\in\mathsf{D}_{\mathsf{per}}(A^\mathsf{e})$ \cite {KS, G1}.

Observe that if $A$ is homologically smooth then so is
$A^{\mathsf{op}}$. Also, if $A$, $B$ are homologically smooth then
their tensor product is a homologically smooth DG
algebra\footnote{There is another important version of the notion of
smoothness in non-commutative geometry, namely, the formal
smoothness (or quasi-freeness) introduced in \cite{CQ}. Notice that
the tensor product of two formally smooth algebras is not a formally
smooth algebra in general.}. In particular, $A^\mathsf{e}$ is
homologically smooth whenever $A$ is. We will use these facts later.

\medskip

\section{Smoothness and saturatedness}

The goal of this section is to prove the following
\begin{theorem}\label{main1}
If $A$ is compact (i.e. $\sum_n
\mathrm{dim}\,\mathrm{H}^n(A)<\infty$) and homologically smooth then
the triangulated category $\mathsf{D}_{\mathsf{per}}(A)$ is
saturated.
\end{theorem}

Let us recall the definition of a saturated triangulated category
\cite{BK1}. Let $\mathcal{T}$ be a triangulated category. It is said
to be Ext-finite if $\sum_n
\mathrm{dim}\,\mathrm{Hom}_{\mathcal{T}}(X,Y[n])<\infty$ for all
$X,Y\in\mathcal{T}$. $\mathcal{T}$ is right (resp. left) saturated
if it is Ext-finite and any contravariant (resp. covariant)
cohomological functor  $h: \mathcal{T}\rightarrow Vect(k)$ of finite
type (i.e. such that $\sum_n\mathrm{dim}\,h(X[n])<\infty$ for any
$X$) is representable. $\mathcal{T}$ is saturated if it is both
right and left saturated.

We will prove the above Theorem for the equivalent category
$\mathsf{Ho}_{\mathsf{per}}(A)$. Let us start with
\begin{lemma}\label{ext}
If $A$ is compact then $\mathsf{Ho}_{\mathsf{per}}(A)$ is Ext-finite.
\end{lemma}
{\bf Proof.} Obviously it is enough to check that $\mathsf{
Ho}(\mathsf{ Tw}(A))$ is Ext-finite.

When $X=A[n]$ and $Y=A[m]$, $\sum_n
\mathrm{dim}\,\mathrm{Hom}_{\mathsf{ Ho}(\mathsf{ Tw}(A))}(X,Y[n])<\infty$
is a straightforward consequence of the compactness of
$A$.  The general case can be reduced to this special case as
follows.

Since $\mathrm{Hom}_{\mathsf{ Ho}(\mathsf{ Tw}(A))}(-, -)$ is exact
with respect to both arguments, it suffices to prove that any object
$X=(\bigoplus\limits_{j=1}^nA[r_j], \alpha)$ of $\mathsf{ Tw}(A)$ is
the cone of a degree 0 morphism between two objects of length at
most $n-1$. Let $\beta=(\alpha_{ij})_{i,j=1,\ldots n-1}$ and
$f=(\alpha_{in})_{i=1,\ldots n-1}$. Then it is easy to see that
$X=\mathrm{Cone}(f)$, where $f$ is viewed as a morphism $(A[r_n-1],
0)\rightarrow(\bigoplus\limits_{j=1}^{n-1}A[r_j], \beta)$. The Lemma is proved.



\medskip
We know that $\mathsf{Ho}_{\mathsf{per}}(A)$ is equivalent to
$\mathsf{Ho}_{\mathsf{per}}(A^{\mathsf{op}})^{\mathsf{op}}$ (see the
proof of Proposition \ref{prop0}). If we can prove that
$\mathsf{Ho}_{\mathsf{per}}(A)$ is right saturated then it would
imply that $\mathsf{Ho}_{\mathsf{per}}(A^{\mathsf{op}})$ is left
saturated. Since we can interchange $A$ and $A^{\mathsf{op}}$ in
this argument ($A^{\mathsf{op}}$ is also compact and homologically
smooth!), we see that to prove Theorem \ref{main1}, it would be
enough to show that $\mathsf{Ho}_{\mathsf{per}}(A)$ is right
saturated.

One has the following result:
\medskip
\begin{theorem}\cite[Theorem 1.3]{BVDB}
Assume $\mathcal{T}$ is an Ext-finite and Karoubian triangulated
category. If it has a strong generator then it is right saturated.
\end{theorem}

Let us explain the last statement. If $\mathcal{T}$ is a
triangulated category and $E\in \mathcal{T}$, define $\langle
E\rangle_n^{\mathcal{T}}$ to be the full subcategory of objects that
can be obtained from $E$ by taking shifts, finite direct sums,
direct summands, and at most $n-1$ cones (details can be found in
\cite[\S 2.1]{BVDB}). An object $E\in\mathcal{T}$ is called a strong
generator if, for some $n$, $\langle E\rangle_n^{\mathcal{T}}$ is
equivalent to $\mathcal{T}$.

The category $\mathsf{Ho}_{\mathsf{per}}(A)$ is Karoubian by its
definition. Let us show that, for a compact and homologically smooth
$A$, $\mathsf{Ho}_{\mathsf{per}}(A)$ is strongly generated.
We are going to use the idea of the proof of Theorem 3.1.4 from \cite[\S 3.4]{BVDB}.

Since $A$ is homologically smooth, there exists a quasi-isomorphism
$pA\rightarrow A$, where
$pA\in\mathsf{Ho}_{\mathsf{per}}(A^\mathsf{e})$ is a direct summand
of some twisted module $(\bigoplus\limits_{j=1}^n A^\mathsf{e}[r_j],
\alpha)$. This means that $pA\in\langle
A^\mathsf{e}\rangle_n^{\mathsf{Ho}_{\mathsf{per}}(A^\mathsf{e})}$.
Take $N\in\mathsf{Ho}_{\mathsf{per}}(A)$. Since $N$ is homotopically
projective, $N\simeq N\dot{\otimes}_AA\simeq N\dot{\otimes}_ApA$
(the latter isomorphism follows from Proposition \ref{flat}).
Therefore
$$
N\simeq N\dot{\otimes}_ApA\in N\dot{\otimes}_A\langle
A^\mathsf{e}\rangle_n^{\mathsf{Ho}_{\mathsf{per}}(A^\mathsf{e})}\subset\langle
N\dot{\otimes}_AA^\mathsf{e}\rangle_n^{\mathsf{Ho}_{\mathsf{per}}(A)}.
$$
(The latter inclusion is due to the fact that the functor
$N\dot{\otimes}_A-:
\mathsf{Ho}_{\mathsf{per}}(A^\mathsf{e})\to\mathsf{Ho}_{\mathsf{per}}(A)$ is triangulated.)
Observe that $N\dot{\otimes}_AA^\mathsf{e}\simeq N\dot{\otimes}_kA$
in $\mathsf{Ho}(A)$. Since each DG $k${\it -module} is homotopically
equivalent to its total cohomology
$\mathrm{H}(N)=\bigoplus\limits_{n\in\mathbb{Z}}\mathrm{H}^n(N)$, we
have $N\dot{\otimes}_k A\simeq\mathrm{H}(N)\dot{\otimes}_kA$ in
$\mathsf{Ho}(A)$.  $\mathrm{H}(N)\dot{\otimes}_kA$ is a free module.
It remains to show that $\mathrm{dim}\,\mathrm{H}(N)<\infty$. The
proof of this fact repeats the proof of Lemma \ref{ext} since the
cohomology functor is cohomological. Theorem \ref{main1} is proved.

\medskip

Let us point out a useful corollary of the above computation:
\begin{proposition}\label{criteria}
Let $A$ be a compact and homologically smooth DG algebra. Then
$N\in\mathsf{D}_{\mathsf{per}}(A)$ iff
$\mathrm{dim}\,\mathrm{H}(N)<\infty$.
\end{proposition}
That $N\in\mathsf{D}_{\mathsf{per}}(A)$ implies
$\mathrm{dim}\,\mathrm{H}(N)<\infty$ is explained in the end of the
above proof. For the converse statement, repeat the above
computation: if $N$ is any module and $pN$ its homotopically
projective resolution then $pN$ is homotopically equivalent to a
direct summand of some $n$-fold extension of the module
$\mathrm{H}(pN)\dot{\otimes}_kA$.

\medskip
\section{Serre duality on $\mathsf{D}_{\mathsf{per}}(A)$ and its applications}\label{hd}

Recall the definition of a Serre functor \cite{BK1}. Let
$\mathcal{T}$ be a $k$-linear Ext-finite triangulated category. A
Serre functor $S: \mathcal{T}\rightarrow\mathcal{T}$ is defined as a
covariant auto-equivalence of $\mathcal{T}$ such that there exists
an isomorphism of bifunctors
\begin{equation}\label{serre}
\mathrm{Hom}_{\mathcal{T}}(X, Y)^*\simeq\mathrm{Hom}_{\mathcal{T}}(Y, S(X)).
\end{equation}
If such a functor exists, it is unique up to an isomorphism.

We want to describe a Serre functor on the category
$\mathsf{D}_{\mathsf{per}}(A)$, where $A$ is compact and
homologically smooth. The answer is known in the case of ordinary
associative algebras (see, for example, \cite [\S 21]{G}); we just
show that the same construction works in our setting. Then we will
compute the inverse functor and prove the main result of the paper,
namely, the existence of a non-degenerate pairing on the Hochschild
homology of a compact homologically smooth DG algebra.

We notice that existence of a Serre functor follows from Theorem
\ref{main1} and the following
\begin{theorem}\cite[\S 3.5]{BK1}
If $\mathcal{T}$ is a saturated triangulated category then it has a Serre functor.
\end{theorem}

Let $S_A$ stand for the so called Nakayama functor on $\mathsf{D}(A)$:
$$
S_{A}: N\mapsto (N^\vee)^*
$$
(see (\ref{dual1}), (\ref{dual2})).
\medskip
\begin{theorem}\label{serreper}
If $A$ is compact and homologically smooth then the functor $S_A$
preserves the subcategory $\mathsf{D}_{\mathsf{per}}(A)\subset
\mathsf{ D}(A)$ and induces a Serre functor on it.
\end{theorem}

Let us explain why $S_A$ preserves $\mathsf{D}_{\mathsf{per}}(A)$.
By Proposition \ref{criteria}, it suffices to prove that $S_A(N)$
has finite dimensional total cohomology whenever $N$ is perfect. As
we know (Proposition \ref{prop0}),
$N^\vee\in\mathsf{D}_{\mathsf{per}}(A^\mathsf{op})$. Then, again by
Proposition \ref{criteria}, the cohomology of the latter module is
finite dimensional, whence the result.

What we are going to show is that $S_A$ is a {\it right} Serre
functor (i.e. it satisfies (\ref{serre}) but is not necessarily an
equivalence). Then Theorem \ref{serreper} will follow from the
existence of a Serre functor on $\mathsf{D}_{\mathsf{per}}(A)$ and
the fact that any two right Serre functors are isomorphic \cite[\S
I.1]{RVDB}.

Let $N,M\in \mathsf{D}_{\mathsf{per}}(A)$.  By (\ref{iso}),
$(\mathsf{RHom}_{A}(N,M))^*\simeq(M\otimes^\mathsf{L}_AN^\vee)^*$.
By (\ref{iso2}),
$(M\otimes^\mathsf{L}_AN^\vee)^*\simeq\mathsf{RHom}_{A}(M,
(N^\vee)^*)$. Theorem \ref{serreper} is proved.

\medskip
It is natural to ask whether the functor $S_A$ can be written in the
form $-\otimes^\mathsf{L}_AX$ for some right DG
$A^\mathsf{e}$-module $X$. To answer this question, consider $A^*$.
It carries a canonical right DG $A^\mathsf{e}$-module structure
coming from the natural left DG $A^\mathsf{e}$-module structure on
$A$:
\begin{equation}\label{lact}
(a'\otimes a'')a=(-1)^{|a'|(|a|+|a''|)}a''aa'.
\end{equation}
Then, by Proposition \ref{prop2}
\medskip
\begin{theorem}
$S_A$ is isomorphic to $-\otimes^\mathsf{L}_AA^*$.
\end{theorem}

\medskip

Let us compute the inverse $S^{-1}_A$. Consider the right DG
$A^\mathsf{e}$-module
$$A^!=\mathsf{RHom}_{(A^\mathsf{e})^{\mathsf{op}}}(A,
(A^\mathsf{e})^{\mathsf{op}}).$$ Here we are using the left DG
$A^\mathsf{e}$-module structure on $A$ defined above (or, more
precisely, the corresponding right DG
$(A^\mathsf{e})^{\mathsf{op}}$-module structure).

\medskip
\begin{theorem}\label{serreinv}
$-\otimes^\mathsf{ L}_AA^!$ is inverse to $S_A$.
\end{theorem}
To show this, it suffices to prove that the two functors form an adjoint pair:
$$
\mathrm{Hom}_{\mathsf{ D}(A)}(N, S_A(M))=\mathrm{Hom}_{\mathsf{
D}(A)}(N\otimes^\mathsf{ L}_AA^!, M),\quad N,M\in
\mathsf{D}_{\mathsf{per}}(A)
$$
which is equivalent to $(\mathrm{Hom}_{\mathsf{D}(A)}(M,
N))^*=\mathrm{Hom}_{\mathsf{ D}(A)}(N\otimes^\mathsf{ L}_AA^!, M)$.
By (\ref{iso2}) and (\ref{twoten}),
$$
\mathsf{RHom}_{A}(N\otimes^\mathsf{L}_AA^!, M)\simeq
(N\otimes^\mathsf{L}_AA^!\otimes^\mathsf{L}_AM^*)^*\simeq
((N\otimes_kM^*)\otimes^\mathsf{L}_{A^\mathsf{e}}A^!)^*.
$$
Proposition \ref{prop1} implies
$((N\otimes_kM^*)\otimes^\mathsf{L}_{A^\mathsf{e}}A^!)^*\simeq\mathsf{RHom}_{A^\mathsf{e}}(A,
N\otimes_kM^*)^*$ and (\ref{homie}) finishes the proof.

\medskip

Thus, we have
\begin{theorem}\label{invert} $A^*$ and $A^!$ are mutually inverse invertible bimodules, i.e. we have isomorphisms
$$
A^*\otimes^\mathsf{L}_{A}A^!\simeq A\simeq
A^!\otimes^\mathsf{L}_{A}A^*
$$
in $\mathsf{D}(A^\mathsf{e})$.
\end{theorem}

\medskip
Now we are ready to prove our main result. Recall \cite{K2} that the
Hochschild homology groups $\mathrm{HH}_n(A)$ are defined as follows
$$
\mathrm{HH}_n(A)=\mathrm{H}^n(A\otimes^\mathsf{L}_{A^\mathsf{e}}A).
$$
The tensor product on the right hand side is defined via the action
(\ref{lact}).

\medskip

\begin{theorem}\label{main} Suppose $A$ is compact and homologically smooth. Then
$\sum_n\mathrm{dim}\,\mathrm{HH}_n(A)<\infty$ and there exists a
canonical non-degenerate pairing
$$
\mathrm{HH}_n(A)\times\mathrm{HH}_{-n}(A)\rightarrow k.
$$
\end{theorem}
Since $A$ is compact and homologically smooth, so is
$(A^\mathsf{e})^\mathsf{op}$. Proposition \ref{criteria} assures
that $A$ is a perfect right DG $(A^\mathsf{e})^\mathsf{op}$-module
(with respect to the action (\ref{lact})). Therefore, by Proposition
\ref{prop0}, there is a canonical isomorphism
$A\simeq\mathsf{RHom}_{A^\mathsf{e}}(A^!, A^\mathsf{e})$. By
Proposition \ref{prop1},
\begin{eqnarray*}
A\otimes^\mathsf{ L}_{A^\mathsf{e}}A\simeq A\otimes^\mathsf{
L}_{A^\mathsf{e}}\mathsf{RHom}_{A^\mathsf{e}}(A^!, A^\mathsf{e})
\simeq \mathsf{RHom}_{A^\mathsf{e}}(A^!, A).
\end{eqnarray*}
In particular, $\sum_n\mathrm{dim}\,\mathrm{HH}_n(A)<\infty$ since
$A^!, A$ are perfect DG $A^\mathsf{e}$-modules (see Lemma
\ref{ext}). Finally, by Theorem \ref{invert} and (\ref{iso2})
$$
\mathsf{RHom}_{A^\mathsf{e}}(A^!,
A)\simeq\mathsf{RHom}_{A^\mathsf{e}}(A^!\otimes^\mathsf{L}_{A}A^*,
A\otimes^\mathsf{L}_{A}A^*)\simeq\mathsf{RHom}_{A^\mathsf{e}}(A,A^*)\simeq(A\otimes^\mathsf{L}_{A^\mathsf{e}}A)^*
$$
and therefore, for any $n$, we have a canonical non-degenerate
pairing
$$\mathrm{H}^n(A\otimes^\mathsf{L}_{A^\mathsf{e}}A)\times\mathrm{H}^{-n}(A\otimes^\mathsf{L}_{A^\mathsf{e}}A)\to
k.$$

\medskip

\appendix
\section{Some canonical isomorphisms}

In this Appendix, we give an account of all canonical isomorphisms used in the paper.

Let $A$ be a DG algebra and $N\in\mathsf{Mod}(A)$,
$M\in\mathsf{Mod}(A^\mathsf{op})$ arbitrary modules. The tensor
product $N\otimes_kM$ is canonically a right DG
$A^\mathsf{e}$-module. Fix a module
$X\in\mathsf{Mod}(A^\mathsf{e})$. One has an obvious isomorphism
\begin{equation}\label{twoten}
N\otimes^\mathsf{L}_AX\otimes^\mathsf{L}_AM\simeq (N\otimes_kM)\otimes^\mathsf{L}_{A^\mathsf{e}}X.
\end{equation}

Let $N$ be a DG $k$-module. Define
\begin{equation}\label{dual1}
N^*=\mathrm{Hom}_{\mathsf{Mod}(k)}(N, k).
\end{equation}
Here $k$ stands for the DG $k$-module whose $0$-th component is $k$
and other components are 0. If $N$ is a right DG module over a DG
algebra then $N^*$ inherits a canonical structure of a right DG
module over the opposite DG algebra. Observe that $N^{**}\simeq N$
on the derived level whenever $N$ has finite dimensional total
cohomology.

Let $N,M\in\mathsf{Mod}(A)$. Both the tensor product $N\otimes_kM^*$
and $A$ are right DG $A^\mathsf{e}$-modules. One has an isomorphism
\begin{equation}\label{homie}
\mathrm{Hom}_{\mathsf{D}(A^\mathsf{e})}(A, N\otimes_kM^*)\simeq\mathrm{Hom}_{\mathsf{D}(A)}(M, N).
\end{equation}

\medskip

Let us now list some canonical isomorphisms involving perfect $A$-modules.
All of them are well known (see, for
example, \cite[\S 21]{G} for a review of the case of ordinary
associative algebras); therefore we give only sketches of proofs.

Take $N\in\mathsf{Mod}(A)$ and set
\begin{equation}\label{dual2}
N^\vee=\mathrm{Hom}_{\mathsf{Mod}(A)}(N, A).
\end{equation}
Thus, $N^\vee\in\mathsf{Mod}(A^\mathsf{ op})$.
\medskip
\begin{proposition}\label{prop0}
The functor $N\mapsto N^\vee$ preserves perfect modules and induces
an equivalence $\mathsf{D}_{\mathsf{per}}(A)\rightarrow
\mathsf{D}_{\mathsf{per}}(A^\mathsf{op})^\mathsf{op}$. More
precisely, there is a canonical isomorphism $N\simeq(N^\vee)^\vee$.
\end{proposition}
Indeed, the functor $N\mapsto\mathrm{Hom}_{\mathsf{Mod}(A)}(N, A)$
is easily seen to preserve twisted modules (the space
$\mathrm{Hom}_{\mathsf{ Tw}(A)}((\oplus_{j=1}^nA[r_j], \alpha), (A,
0))$ consists of $n\times1$-matrices of elements of $A$($=A^\mathsf{
op}$), and the differential and the $A$-action are exactly of the
same form as those in $\mathsf{ Tw}(A^\mathsf{ op})$). It clearly
descents to a functor from $\mathsf{Ho}_{\mathsf{per}}(A)$ to
$\mathsf{Ho}_{\mathsf{per}}(A^{\mathsf{op}})^{\mathsf{op}}$ since
$\mathsf{Hom}_A(-, A)$ sends direct summands to direct summands.
Furthermore, one has an obvious canonical map
$$
N\rightarrow\mathsf{Hom}_{A^{\mathsf{op}}}(\mathsf{Hom}_{A}(N, A), A^\mathsf{op}).
$$
It is easy to see that this map is an isomorphism when $N$ is a
twisted module. If $N$ is a direct summand of a twisted module $P$,
we have $$ P\simeq\mathsf{Hom}_{A^{\mathsf{op}}}(\mathsf{Hom}_{A}(P,
A), A^\mathsf{op}), \quad
N\hookrightarrow\mathsf{Hom}_{A^{\mathsf{op}}}(\mathsf{Hom}_{A}(N,
A), A^\mathsf{op})
$$
which implies
$N\simeq\mathsf{Hom}_{A^{\mathsf{op}}}(\mathsf{Hom}_{A}(N, A),
A^\mathsf{op})$. To finish the proof, it remains to observe that,
for an arbitrary perfect module $N$,
$(N^\vee)^\vee\simeq\mathsf{Hom}_{A^{\mathsf{op}}}(\mathsf{Hom}_{A}(pN,
A), A^\mathsf{op})$, where $pN$ is a
resolution of $N$.

\medskip
\begin{proposition}\label{prop1}
If $N$ is perfect and $M$ is arbitrary then there is a canonical isomorphism
\begin{equation}\label{iso}
\mathsf{RHom}_{A}(N,M)\simeq M\otimes^\mathsf{ L}_A N^\vee.
\end{equation}
\end{proposition}
To prove this Proposition, consider the map of DG modules
$$
F_{N,M}: M{\otimes}_A \mathrm{Hom}_{\mathsf{ Mod}(A)}(N,
A)\rightarrow \mathrm{Hom}_{\mathsf{ Mod}(A)}(N,M)
$$
given by $F_{N,M}(m\otimes f)(n)=mf(n)$, where $m\otimes f\in
M{\otimes}_A \mathrm{Hom}_{\mathsf{ Mod}(A)}(N, A)$ and $n\in N$. It
induces a map
$$
{\bar{F}}_{N,M}: M\dot{\otimes}_A \mathsf{Hom}_{A}(N, A)\rightarrow \mathsf{Hom}_{A}(N,M).
$$
If $N$ is a twisted module then a straightforward computation shows
that $F_{N,M}$ is bijective, and so is ${\bar{F}}_{N,M}$. If $N$ is
a direct summand of a twisted module $P$ then one can show that the
image of $M\dot{\otimes}_A \mathsf{Hom}_{A}(N, A)$ under
${\bar{F}}_{P,M}$ is inside of
$\mathsf{Hom}_{A}(N,M)\subset\mathsf{Hom}_{A}(P,M)$.

\medskip

Finally, we want to mention the following
\begin{proposition}\label{prop2}
Suppose $N$ is a perfect right DG $A$-module and $M$ is an arbitrary
right DG $A\otimes B$-module, where $B$ is yet another DG algebra.
Then there is a natural isomorphism of right DG
$B^{\mathsf{op}}$-modules
\begin{equation}\label{iso2}
(\mathsf{RHom}_{A}(N,M))^*\simeq N\otimes^\mathsf{ L}_A M^*.
\end{equation}
\end{proposition}
The proof of this proposition is completely analogous to the proof
of the preceding one. This time the isomorphism comes from the
canonical map
$$
G_{N,M}: N\otimes_A M^*\rightarrow (\mathrm{Hom}_{\mathsf{Mod}(A)}(N,M))^*
$$
given by $G_{N,M}(n\otimes \nu)(f)=(-1)^{|n|(|\nu|+|f|)}\nu(f(n))$,
where $n\otimes \nu\in N\otimes_A M^*$ and
$f\in\mathrm{Hom}_{\mathsf{Mod}(A)}(N,M)$. If $N$ is a direct
summand of a twisted $A$-module then $G_{N,M}$ is
obviously bijective.

\small{\it Mathematics Department, Kansas State University}

\small{\it 138 Cardwell Hall}

\small{\it Manhattan, KS 66506-2602}

\small{\it e-mail: shklyarov@math.ksu.edu}

\medskip




\medskip
\begin{thebibliography}{99}
\bibitem{BL} J. Bernstein, V. Lunts, Equivariant sheaves and functors. Lecture Notes in Mathematics, 1578. Springer-Verlag, Berlin, 1994.
\bibitem{BK}  A. Bondal,  M. Kapranov, Enhanced triangulated categories. Math. USSR-Sb.  70  (1991),  no. 1, 93--107.
\bibitem{BK1} A. Bondal,  M. Kapranov,  Representable functors, Serre functors, and reconstructions. Math. USSR-Izv. 35 (1990), no. 3, 519--541.
\bibitem{BVDB} A. Bondal, M. Van den Bergh, Generators and representability of functors in commutative and noncommutative geometry. Preprint: math.AG/0204218.
\bibitem{Ca} A. Caldararu, The Mukai pairing, I: the Hochschild
structure. Preprint: math.AG/0308079.
\bibitem{Ci} C. Cibils, Hochschild homology of an algebra whose quiver has
no oriented cycles.  Representation theory, I (Ottawa, Ont., 1984),
55--59, Lecture Notes in Math., 1177, Springer, Berlin, 1986.
\bibitem{C} K. Costello, Topological conformal field theories and Calabi-Yau
categories. Preprint: math.QA/0412149.
\bibitem{CQ} J. Cuntz, D. Quillen, Algebra extensions and nonsingularity.  J. Amer. Math. Soc.  8  (1995),  no. 2, 251--289.
\bibitem{GJ} E. Getzler, J. D. S. Jones, $A_\infty$-algebras and the cyclic bar complex. Illinois J. Math. 34 (1990), no. 2, 256--283.
\bibitem{G} V. Ginzburg, Lectures on noncommutative geometry. Preprint: math.AG/0506603.
\bibitem{G1} V. Ginzburg, Calabi-Yau algebras. Preprint:
math.AG/0612139.
\bibitem{D} V. Drinfeld, DG quotients of DG categories. Preprint: math.KT/0210114.
\bibitem{Ka} D. Kaledin, Non-commutative Hodge-to-de Rham degeneration via the method of
Deligne-Illusie. Preprint: math.KT/0611623.
\bibitem{K} B. Keller, Deriving DG categories. Ann. Sci. \'{E}cole Norm. Sup. (4) 27 (1994), no. 1, 63--102.
\bibitem{K1} B. Keller On the cyclic homology of ringed spaces and schemes.  Doc. Math.  3  (1998), 231--259.
\bibitem{K11/2} B. Keller, On the cyclic homology of exact categories.  J. Pure Appl. Algebra  136  (1999),  no. 1, 1--56.
\bibitem{K2} B. Keller, Invariance and localization for cyclic homology of DG algebras.  J. Pure Appl. Algebra  123  (1998),  no. 1-3, 223--273.
\bibitem{K3} B. Keller, On differential graded categories. Preprint:
math.KT/0601185.
\bibitem{KS} M. Kontsevich, Y. Soibelman, Notes on A-infinity algebras, A-infinity categories and non-commutative geometry. I Preprint: math.RA/0606241.
\bibitem{Ma0} N. Markarian, Poincare-Birkhoff-Witt isomorphism, Hochschild homology and Riemann-Roch theorem.
Preprint: MPI 2001-52.
\bibitem{Ma} N. Markarian, The Atiyah class, Hochschild cohomology and
the Riemann-Roch theorem. Preprint: math.AG/0610553.
\bibitem{RVDB} I. Reiten, M. Van den Bergh, Noetherian hereditary abelian categories satisfying Serre duality. math.RT/9911242.
\bibitem{Ro} R. Rouquier, Dimensions of triangulated categories.
Preprint: math.CT/0310134.
\bibitem{TV} B. Toen, M. Vaquie, Moduli of objects in dg-categories.
Preprint: math.AG/0503269.
\bibitem{vdb} M. van den Bergh, A relation between Hochschild homology and
cohomology for Gorenstein rings.  Proc. Amer. Math. Soc.  126
(1998),  no. 5, 1345--1348.
\end{thebibliography}
\end{document}